\pdfoutput=1
\documentclass[12pt]{amsart}

\usepackage{float}
\usepackage{amsmath,amsfonts,amssymb,mathabx,latexsym,mathtools}
\usepackage{enumitem}
\usepackage[usenames, dvipsnames]{xcolor}
\usepackage{hyperref}
\usepackage{ytableau}
\usepackage{centernot}
\usepackage{tikz-cd}
\usepackage{xspace}
\usepackage{tikz}
\usetikzlibrary{decorations.markings}
\usepackage{stmaryrd}
\usepackage{wrapfig}
 \usepackage{url}

\usetikzlibrary{calc, shapes, backgrounds,arrows,positioning,plotmarks}
\tikzset{>=stealth',
  head/.style = {fill = white, text=black},
  plaque/.style = {draw, rectangle, minimum size = 10mm}, 
  pil/.style={->,thick},
  junct/.style = {draw,circle,inner sep=0.5pt,outer sep=0pt, fill=black}
  }

\newtheorem{theorem}{Theorem}[section]

\newtheorem{claim}[theorem]{Claim}
\newtheorem{proposition}[theorem]{Proposition}

\newtheorem{conjecture}[theorem]{Conjecture}

\theoremstyle{definition}

\newtheorem{examplex}[theorem]{Example}
\newenvironment{example}
  {\pushQED{\qed}\examplex}
  {\popQED\endexamplex}

\theoremstyle{remark}
\newtheorem{remark}[theorem]{Remark}

\numberwithin{equation}{section}

\newcommand{\ghost}{{\includegraphics[height=1.4ex]{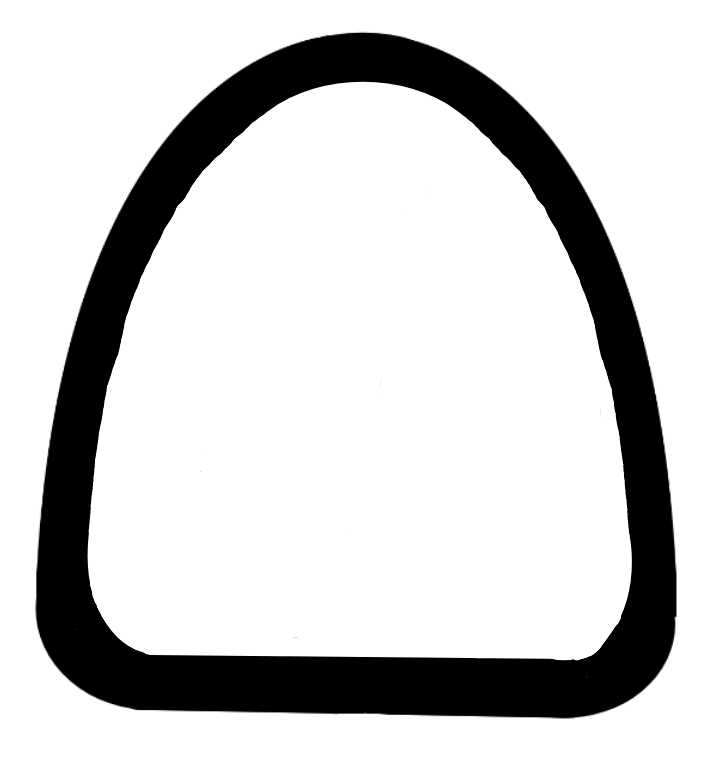}}}

\newcommand\mydef[1]{{\bf #1}}
\newcommand\wt{{\tt wt}}

\newcommand\Koh[1]{{\sf Koh}(#1)}
\newcommand\KKoh[1]{{{\sf KKoh}}(#1)}
\newcommand\bKKoh[1]{{{\sf boldKKoh}}(#1)}
\newcommand\tKKoh[1]{\overline{{\sf KKoh}}(#1)}
\newcommand\al[1]{\alpha_{#1}}
\newcommand\be[1]{\beta_{#1}}
\newcommand\ob{{\bullet}}
\newcommand\gb{{\ghost}}
\begin{document}

\title[A counterexample to Ross--Yong conjecture for Grothendieck polynomials]{A counterexample to the Ross--Yong conjecture for Grothendieck polynomials}

\author{Colleen Robichaux}
\address[CR]{
Dept.~of Mathematics,
University of California, Los Angeles,
Los Angeles, CA  \newline \indent 90095, USA}
\email{robichaux@math.ucla.edu}
\thanks{Colleen Robichaux was supported by the NSF MSPRF No. DMS 2302279.}

\date{\today}

\begin{abstract}
We give a minimal counterexample for a conjecture of Ross and Yong (2015) which proposes a K-Kohnert rule for Grothendieck polynomials. We conjecture a revised version of this rule. We then prove both rules hold in the $321$-avoiding case.
\end{abstract}

\maketitle

\section{Introduction}
Introduced by A.~Lascoux and M.~P.~Sch\"utzenberger \cite{LS82} to study of the K-theory of the complete flag variety,
\emph{Grothendieck polynomials} ${\mathfrak G}_{w}$ are a basis for the polynomial ring $\mathbb{Z}[x_1,x_2, \ldots]$. Additionally A.~Lascoux and M.~P.~Sch\"utzenberger \cite{LS1}, introduce \emph{Schubert polynomials} ${\mathfrak S}_{w}$, another basis for $\mathbb{Z}[x_1,x_2, \ldots]$,  comprised of the lowest degree terms of ${\mathfrak G}_{w}$.
Several combinatorial rules have been developed to study Grothendieck polynomials \cite{BFHTW,BS,FominKrillov,ice,BPD}. Many of these naturally generalize formulas for Schubert polynomials. 

Kohnert's rule \cite{AABE,Assaf,Kohnert, Winkel1,Winkel2} combinatorially computes Schubert polynomials through local moves on diagrams.
C.~Ross and A.~Yong \cite{Ghost} conjecture a generalized Kohnert's rule  to compute Grothendieck polynomials and \emph{Lascoux polynomials}, another basis for $\mathbb{Z}[x_1,x_2, \ldots]$.
 After initial work of O.~Pechenik and T.~Scrimshaw \cite{PS},
 J.~Pan and T.~Yu \cite{LascKKoh} prove the Ross--Yong rule for Lascoux polynomials. 
 
 In contrast to the Lascoux case, 
no progress has been made towards establishing the Ross--Yong conjecture for Grothendieck polynomials. 
Diagrams in recent work studying properties of Grothendieck polynomials \cite{bubbleGroth,MSS,pan2023topdegree} are noted as similar to those in the conjecture. The Ross--Yong conjecture has been checked through $S_7$, beyond which a counterexample might be surprising. 

In Section~\ref{sec:RY} we present a minimal counterexample in $S_8$ to the Ross--Yong conjecture for Grothendieck polynomials. This counterexample suggests that when studying Grothendieck polynomials, behavior exhibited in small examples may be misleading.
In Section~\ref{sec:newRule} we propose a new K-theoretic Kohnert rule for Grothendieck polynomials. This new rule has been checked through $S_8$ and more closely resembles the pipe dream rule for Grothendieck polynomials \cite{FominKrillov}. 
We end with a proof for the $321$-avoiding case of both rules in Section~\ref{sec:321pf}.

\section{Background}
Let $S_n$ denote the \mydef{symmetric group} on $n$ elements.
We say $w\in S_n$ contains a pattern $p$ if some subsequence in $w$ has the same relative order as the entries in $p$. For example, $w=1\underline{7}\underline{4}6\underline{2}35$ contains a $321$ pattern, which we have underlined.
If $w$ does not contain a pattern $p$, we say $w$ is \mydef{$p$-avoiding}.

For $f\in\mathbb{Z}[x_1,x_2,\ldots,x_n]$, define
 \[\partial_if=\frac{f-s_if}{x_i-x_{i+1}} \text{ and } \widetilde{\pi}_if=\partial_i(1-x_{i+1})f.\]
Here $s_i$ is the transposition swapping $i$ and $i+1$. 
For $n\in\mathbb{Z}$, let $[n]:=\{i\in\mathbb{Z}  \, | \, 1\leq i\leq n\}$.

Indexed by $w\in S_n$,
\mydef{Schubert polynomials} ${\mathfrak S}_{w}(x_1,\ldots,x_n)$ and \mydef{Grothendieck polynomials} ${\mathfrak G}_{w}(x_1,\ldots,x_n)$ are defined recursively \cite{LS1, LS82}. 
For $w_0=n \ n-1 \ \cdots 2 \ 1 $,
\begin{align*}
     {\mathfrak S}_{w_0}(x_1,\ldots,x_n)={\mathfrak G}_{w_0}(x_1,\ldots,x_n)&:=x_1^{n-1}x_2^{n-2}\cdots x_{n-1}.
\end{align*}
Otherwise for $w\neq w_0$, take $i\in[n]$ such that $w(i)<w(i+1)$. Then 
\begin{align*}
     \mathfrak S_{w}(x_1,\ldots,x_n)&=\partial_i\mathfrak S_{ws_i}(x_1,\ldots,x_n), \text{ and }\\
     \mathfrak G_{w}(x_1,\ldots,x_n)&=\widetilde{\pi}_i\mathfrak G_{ws_i}(x_1,\ldots,x_n).
\end{align*}

The \mydef{Lascoux polynomials} ${\mathfrak L}_{\alpha}(x_1,\ldots,x_n)$ are defined recursively for weak compositions $\alpha\in\mathbb{Z}_{\geq 0}^n$. 
If $\alpha=(\alpha_1,\ldots,\alpha_n)$ is weakly decreasing, \[{\mathfrak L}_{\alpha}(x_1,\ldots,x_n)=x^{\alpha}=x_1^{\alpha_1} \cdots x_n^{\alpha_n}.\] 
Otherwise, take $i\in[n]$ such that $\alpha_i<\alpha_{i+1}$. Then
\[{\mathfrak L}_{\alpha}(x_1,\ldots,x_n)=\widetilde{\pi}_i((1-x_{i+1}){\mathfrak L}_{s_i\alpha}(x_1,\ldots,x_n)),\] 
where $s_i\alpha=(\alpha_1,\ldots,\alpha_{i-1},\alpha_{i+1},\alpha_{i},\ldots,\alpha_n)$.

\subsection{Pipe Dreams} 
For $P\subseteq [n]\times[n]$ define the \mydef{weight} of $P$ as
 $\wt(P)\in\mathbb{Z}_{\geq0}^n$ such that $\wt(P)_i:=\# \{(i,c)\in P \, | \,  c\in[n]\}$, where $i\in[n]$.
Then take $\#P=\sum_{i\in[n]}\wt(P)_i$.

Label boxes of $[n]\times[n]$ where $(i,k)$ has label $i+k-1$.
 For $P\subseteq [n]\times[n]$, let ${\sf word}(P)$ be the sequence determined by recording the labels of $P$, reading right to left across rows in $[n]\times[n]$, starting with the top row and working downwards. For example, if $P=[3]\times[3]$, then ${\sf word}(P)=321432543$.
 Note $(1,1)$ corresponds to the northwestmost box in $[n]\times [n]$. 
 
 Define an algebra over $\mathbb Z$ generated by $\{e_w \, | \, w\in S_n\}$, where multiplication is given by
\[e_w\cdot e_{s_i}=
\begin{cases}
e_{ws_i} &\text{if } \ell(ws_i)>\ell(w), \\
e_w & \text{otherwise}.
\end{cases}\]
 The \mydef{Demazure product} $\delta(P)$ of $P$ is the permutation that satisfies \[e_{s_{i_1}}\cdots e_{s_{i_j}}=e_{\delta(P)},\]
 where ${\sf word}(P)=(i_1,i_2,\ldots,i_j)$. 
 Following the setup in \cite{KM}, define
 \begin{align*} 
     {\sf{Pipes}}(w)&:=\{P\subseteq [n]\times[n]\,|  \,\delta(P)=w\}.\\
\end{align*}
  We illustrate $P\in {\sf{Pipes}}(w)$ by placing a $+$ at each $(i,j)\in P$ in the $[n]\times[n]$ grid. We call ${\sf{Pipes}}(w)$ the set of \mydef{pipe dreams} for $w$. Minimal pipe dreams generate Schubert polynomials. 
\begin{theorem}\cite{BJS,FKSch}\label{thm:rPD}
    For $w\in S_n$,
     \[{\mathfrak{S}_w}= \sum_{\substack{P\in {\sf{Pipes}}(w) \\ \#P=\ell(w)}}x^{\wt(P)}.\]
\end{theorem}

  This Schubert formula naturally generalizes to Grothendieck polynomials.
\begin{theorem}\cite{FominKrillov}\label{thm:PD}
    For $w\in S_n$,
     \[{\mathfrak{G}_w}= \sum_{P\in {\sf{Pipes}}(w)}(-1)^{\#P-\ell(w)}x^{\wt(P)}.\]
\end{theorem}

Notice ${\mathfrak{G}_w}=\sum_{\gamma}(-1)^{|\gamma|-\ell(w)}g_{w,\gamma}x^{\gamma}$, where
\[g_{w,\gamma}=\#\{P\in {\sf{Pipes}}(w) \, | \, \wt(P)=\gamma\}.\]

\begin{example}\label{ex:PD}
    Take $w=12365847$ and $\gamma=(3,3,3,2)$. Below are three $P\in {\sf{Pipes}}(w)$ such that $ \wt(P)=\gamma$. The words associated to each are $75475475464,75475476464,75475475476$, respectively. Thus, $g_{w,\gamma}\geq 3$.
      \[\begin{picture}(300,90)
\put(0,75){\ytableausetup
{boxsize=0.8em}
{\begin{ytableau}
 \ &  \ &  \ &  + &  + & \ &  + &  \ \\
 \ &  \ &  + &  + & \ &  + &  \ &  \ \\
 \ &  + &  + & \ &  + & \ &  \ &  \ \\
 + &  \ &  + &  \ &  \ & \ &  \ &  \ \\
 \ &  \ &  \ &  \ &  \ & \ &  \ &  \ \\
 \ &  \ &  \ &  \ &  \ & \ &  \ &  \ \\
 \ &  \ &  \ &  \ &  \ & \ &  \ &  \ \\
 \ &  \ &  \ &  \ &  \ & \ &  \ &  \ \\
\end{ytableau}}}
\put(100,75){\ytableausetup
{boxsize=0.8em}
{\begin{ytableau}
 \ &  \ &  \ &  + &  + & \ &  + &  \ \\
 \ &  \ &  + &  + & \ &  + &  \ &  \ \\
 \ &  + &  \ & + &  + & \ &  \ &  \ \\
 + &  \ &  + &  \ &  \ & \ &  \ &  \ \\
 \ &  \ &  \ &  \ &  \ & \ &  \ &  \ \\
 \ &  \ &  \ &  \ &  \ & \ &  \ &  \ \\
 \ &  \ &  \ &  \ &  \ & \ &  \ &  \ \\
 \ &  \ &  \ &  \ &  \ & \ &  \ &  \ \\
\end{ytableau}}}
\put(200,75){\ytableausetup
{boxsize=0.8em}
{\begin{ytableau}
 \ &  \ &  \ &  + &  + & \ &  + &  \ \\
 \ &  \ &  + &  + & \ &  + &  \ &  \ \\
 \ &  + &  + & \ &  + & \ &  \ &  \ \\
 \ &  \ &  + &  + &  \ & \ &  \ &  \ \\
 \ &  \ &  \ &  \ &  \ & \ &  \ &  \ \\
 \ &  \ &  \ &  \ &  \ & \ &  \ &  \ \\
 \ &  \ &  \ &  \ &  \ & \ &  \ &  \ \\
 \ &  \ &  \ &  \ &  \ & \ &  \ &  \ \\
\end{ytableau}}}
\end{picture}
\]
 In fact,
    these are all the $P\in {\sf{Pipes}}(w)$ with weight $\gamma$, so $g_{w,\gamma}=3$.
\end{example}

\subsection{Kohnert Diagrams}
A \mydef{diagram} is any $D\subseteq [n]\times [n]$. Here $(1,1)$ corresponds to the northwestmost box in $[n]\times [n]$. Let $\#D$ denote the cardinality of the set $D$. 
We illustrate $D$ by marking boxes of $D$ in $[n]\times [n]$ with $\ob$. 
The \mydef{weight} of $D$, denoted $\wt(D)\in \mathbb{Z}_{\geq 0}^n$, is defined by
\[\wt(D)_i:=\# \text{ of boxes in row } i, \text{ where } i\in[n].\]

For a diagram $D$ and $(i,j)\in D$ rightmost in some row $i\in[n]$, the \mydef{Kohnert move} on $D$ at $(i,j)$ outputs $D'=D-\{(i,j)\}\cup \{(i',j)\}$, where $i'=\max\{r\in[i] \, | \, (r,j)\not\in D\}$. 
Let $\Koh{D}$ denote the set of all diagrams obtainable through applying successive Kohnert moves on $D$.

The \mydef{Rothe diagram} of $w\in S_n$ is the set
	\[ D(w):= \{(i, j)\in [n]\times [n] \,|  \,  w_i > j \text{ and } w^{-1}_j > i\}.\]
 The \mydef{key diagram} of $\alpha\in \mathbb{Z}_{\geq0}^n$ is the set
	\[ D(\alpha):= \{(i, j)\in [n]\times [n] \,|  \,  \alpha_i \geq j\}.\]

 A diagram $D$ is \mydef{northwest hook-closed} if $(i,j),(i',j')\in D$ such that $i'<i, j'>j$ implies $(i',j)\in D$. A diagram $D$ is \mydef{southwest hook-closed} if $(i,j),(i',j')\in D$ such that $i'>i, j'>j$ implies $(i',j)\in D$. Both Rothe diagrams and key diagrams are northwest hook-closed. Key diagrams are southwest hook-closed, but Rothe diagrams may not be. 

\begin{example}\label{ex:Rothe} Take $w=418723956$. On the left is $D(w)$. On the right is some diagram $D\in \Koh{D(w)}$.
\[\begin{picture}(250,60)
\put(0,55){\ytableausetup
{boxsize=0.8em}
{\begin{ytableau}
 \ob &  \ob &  \ob &  \ &  \ & \ &  \ &  \ &  \ \\
 \ &  \ &  \ &  \ &  \ & \ &  \ &  \ &  \ \\
 \ &  \ob &  \ob &  \ &  \ob & \ob &  \ob &  \ &  \ \\
 \ &  \ob &  \ob &  \ &  \ob & \ob &  \ &  \ &  \ \\
\ &  \ &  \ &  \ &  \ & \ &  \ &  \ &  \ \\
\ &  \ &  \ &  \ &  \ & \  &  \ &  \ & \  \\
  \ &  \ & \ &  \ &  \ob & \ob &  \ &  \ &  \ \\
 \ &  \ &  \ &  \ &  \ & \ &  \ &  \ &  \ \\
  \ &  \ &  \ &  \ &  \ & \ &  \ &  \ &  \ \\
\end{ytableau}}}
\put(150,55){\ytableausetup
{boxsize=0.8em}
{\begin{ytableau}
 \ob &  \ob &  \ob &  \ &  \ob & \ob &  \ &  \ &  \ \\
 \ &  \ &  \ &  \ &  \ob & \ob &  \ob &  \ &  \ \\
 \ &  \ob &  \ob &  \ &  \ob & \ob &  \ &  \ &  \ \\
 \ &  \ob &  \ob &  \ &  \ & \ &  \ &  \ &  \ \\
\ &  \ &  \ &  \ &  \ & \ &  \ &  \ &  \ \\
\ &  \ &  \ &  \ &  \ & \  &  \ &  \ & \  \\
  \ &  \ & \ &  \ &  \ & \ &  \ &  \ &  \ \\
 \ &  \ &  \ &  \ &  \ & \ &  \ &  \ &  \ \\
  \ &  \ &  \ &  \ &  \ & \ &  \ &  \ &  \ \\
\end{ytableau}}}
\end{picture}
\]
\end{example}

These Kohnert diagrams give another rule to compute ${\mathfrak{S}_w}$. 
\begin{theorem}\cite{Assaf,Kohnert}\label{thm:Kohnert}
 For $w\in S_n$,
    \[{\mathfrak{S}_w}= \sum_{D\in \Koh{D(w)}}x^{\wt(D)}.\]
\end{theorem}
See \cite{Lenart} for discussion of early efforts to prove Theorem~\ref{thm:Kohnert}. This includes a conjectural bijection posed by S.~Billey between $\Koh{D(w)}$ and $\{P\in {\sf{Pipes}}(w)  \, | \, \#P=\ell(w)\}$. More recent proofs have been proposed by S.~Assaf \cite{Assaf} and S.~Armon, S.~Assaf, G.~Bowling, and H.~Ehrhard \cite{AABE}.

\section{The Ross--Yong conjecture}\label{sec:RY}
Let $D\subseteq[n]\times [n]$ be a diagram. In this setting, boxes in $D$ are labeled either $\ob$ or $\gb$. For ease of reading we may omit southmost empty rows and eastmost empty columns in the grid.
Suppose $(i,j)\in D$ is rightmost in some row $i\in[n]$ and has label $\ob$. Take \[i'=\max\{r\in[i] \, | \, (r,j)\not\in D\}.\] Then if each $(r,j)\in D$ for $i'+1\leq r \leq i$ has label $\ob$, define the following:
\begin{itemize}
    \item[--] The \mydef{Kohnert move} on $D$ at $(i,j)$ outputs the diagram $D'=D-\{(i,j)\}\cup \{(i',j)\}$. The new box $(i',j)$ has label $\ob$.
    \item[--] The \mydef{K-Kohnert move} on $D$ at $(i,j)$ outputs the diagram $D'=D\cup \{(i',j)\}$. The new box $(i',j)$ has label $\ob$, and the box $(i,j)$ is reassigned label $\gb$.
\end{itemize}
\begin{example}
    Let $D$ be the diagram to the left. The middle diagram is the output of applying a Kohnert move to $(5,1)$. The rightmost diagram is the output of applying a K-Kohnert move to $(5,1)$.
\[\begin{picture}(250,45)
\put(0,30){\ytableausetup
{boxsize=1em}
{\begin{ytableau}
 \ &  \  &  \ \\
      \ &  \  &  \ \\
       \ob &  \ob  &  \ \\
      \ob &  \ & \ob \\
      \ob &  \  & \  \\
\end{ytableau}}}
\put(100,30){\ytableausetup
{boxsize=1em}
{\begin{ytableau}
 \ &  \  &  \ \\
     \ob  &  \  &  \ \\
     \ob &  \ob  &  \ \\
      \ob &  \ & \ob \\
      \ &  \  & \  \\
\end{ytableau}}}
\put(200,30){\ytableausetup
{boxsize=1em}
{\begin{ytableau}
 \ &  \  &  \ \\
     \ob &  \  &  \ \\
     \ob &  \ob  &  \ \\
      \ob &  \ & \ob \\
      \gb &  \  & \  \\
\end{ytableau}}}
\end{picture}
\]
\end{example}

Let $\KKoh{D}$ denote the set of all diagrams obtainable through successive Kohnert moves and K-Kohnert moves on $D$, where the initial diagram $D$ has all boxes labelled $\ob$. C.~Ross and A.~Yong \cite{Ghost} conjectured these diagrams generate Lascoux and Grothendieck polynomials.\footnote{We note that \cite[Conjecture 1.6]{Ghost} was misstated in the 2015 version. The statement of this conjecture was updated in 2018. This 2018 statement is consistent with C.~Ross's 2011 report \cite{RY11} that originally stated their rule.}

 J.~Pan and T.~Yu \cite{LascKKoh} prove this conjecture through a bijection to particular set-valued tableaux defined by M.~Shimozono and T.~Yu \cite{SU}.
\begin{theorem}\cite{LascKKoh}\label{thm:Lasc} For $\alpha,\gamma\in\mathbb{Z}_{\geq 0}^n$, \[[x^\gamma]{\mathfrak{L}_\alpha}=(-1)^{|\gamma|-|\alpha|}\#\{D\in \KKoh{D(\alpha)} \, | \, \wt(D)=\gamma\}.\] 
\end{theorem}

Recall that $g_{w,\gamma}=\#\{P\in {\sf{Pipes}}(w) \, | \, \wt(P)=\gamma\}$. Then define
 \[\KKoh{w,\gamma}:=\{D\in \KKoh{D(w)} \, | \, \wt(D)=\gamma\}.\] 

\begin{conjecture}\cite[Conjecture 1.6]{Ghost}\label{conj:main} For $w\in S_n$ and $\gamma\in\mathbb{Z}_{\geq 0}^n$, 
\[g_{w,\gamma}=\#\KKoh{w,\gamma}.\]
\end{conjecture}

As stated in \cite{Ghost}, Conjecture~\ref{conj:main} holds for $n\leq 7$. This conjecture fails in $S_8$.
\begin{claim}\label{claim:main}
    Conjecture~\ref{conj:main} is false.
\end{claim}
\begin{proof}
     Take $w=12365847$.  As established in Example~\ref{ex:PD}, $g_{w,\gamma}=3$ for $\gamma=(3,3,3,2)$. 
     We prove the following:
     \begin{equation}\label{eq:failConj}
        \#\KKoh{w,\gamma}=2<g_{w,\gamma}.
     \end{equation}
This counterexample has been confirmed computationally. For completeness, we include a non-computer based proof here.

To the left is $D(w)$. To the right are diagrams $D_1,D_2\in \KKoh{w,\gamma}$.
\begin{equation*}
\begin{picture}(350,90)
\put(0,70){\ytableausetup
{boxsize=1em}
{\begin{ytableau}
 \ &  \ &  \ &  \ &  \ & \ &  \ &  \ \\
   \ &  \ &  \ &  \ &  \ & \ &  \ &  \ \\
    \ &  \ &  \ &  \ &  \ & \ &  \ &  \ \\
     \ &  \ &  \ &  \ob &  \ob & \ &  \ &  \ \\
   \ &  \ &  \ &  \ob &  \ & \ &  \ &  \ \\
    \ &  \ &  \ &  \ob &  \ & \ &  \ob &  \ \\
\end{ytableau}}}
\put(125,70){\ytableausetup
{boxsize=1em}
{\begin{ytableau}
 \ &  \ &  \ &  \ob &  \ob & \ &  \ob &  \ \\
   \ &  \ &  \ &  \gb &  \gb & \ &  \gb &  \ \\
    \ &  \ &  \ &  \ob &  \gb & \ &  \gb &  \ \\
     \ &  \ &  \ &  \ob &  \gb & \ &  \ &  \ \\
   \ &  \ &  \ &  \ &  \ & \ &  \ &  \ \\
    \ &  \ &  \ &  \ &  \ & \ &  \ &  \ \\
\end{ytableau}}}
\put(250,70){\ytableausetup
{boxsize=1em}
{\begin{ytableau}
 \ &  \ &  \ &  \ob &  \ob & \ &  \ob &  \ \\
   \ &  \ &  \ &  \gb &  \gb & \ &  \gb &  \ \\
    \ &  \ &  \ &  \ob &  \gb & \ &  \gb &  \ \\
     \ &  \ &  \ &  \ob &  \ & \ &  \gb &  \ \\
   \ &  \ &  \ &  \ &  \ & \ &  \ &  \ \\
    \ &  \ &  \ &  \ &  \ & \ &  \ &  \ \\
\end{ytableau}}}
\end{picture}
\end{equation*}
Suppose $D\in \KKoh{w,\gamma}-\{D_1,D_2\}$. The number of boxes in each column of $D(w)$ as well as the condition that $\ob$'s cannot pass over $\gb$'s
ensures $U\subset D$, where $U$ is the diagram below. 
In $U$ we have included row and column indices for ease of reading: 
\[\begin{picture}(100,95)
\put(5,70){\ytableausetup
{boxsize=1em}
{\begin{ytableau}
\ &  \ &  \ &  \ob &  \ob & \ &  \ob &  \ \\
   \ &  \ &  \ &  \ &  \gb & \ &  \gb &  \ \\
    \ &  \ &  \ &  \ &  \gb & \ &  \gb &  \ \\
     \ &  \ &  \ &  \ & \ & \ &  \ &  \ \\
   \ &  \ &  \ &  \ &  \ & \ &  \ &  \ \\
    \ &  \ &  \ &  \ &  \ & \ &  \ &  \ \\
\end{ytableau}}}
\put(0,72){\small$ 1$}
\put(0,60){\small$ 2$}
\put(0,48){\small$ 3$}
\put(0,35){\small$ 4$}
\put(0,23){\small$ 5$}
\put(0,11){\small$ 6$}
\put(12,85){\small$ 1$}
\put(24,85){\small$ 2$}
\put(36,85){\small$ 3$}
\put(48,85){\small$ 4$}
\put(61,85){\small$ 5$}
\put(73.5,85){\small$ 6$}
\put(86,85){\small$ 7$}
\put(98,85){\small$ 8$}
\end{picture} \ \ \ \ \raisebox{0.5cm}{.}
\]

Since $\wt(D)=\gamma$, we have $(2,4),(3,4)\in D$ and 
\begin{equation}\label{eq:fillBox}
    \#(D\cap \{(4,4),(4,5),(4,7)\})=2.
\end{equation}

 If $(2,4)$ has label $\gb$ in $D$, then $(3,4),(4,4)\in D$ must have label $\ob$. However this implies $D\in\{D_1,D_2\}$, a contradiction. 

Now assume $(2,4)$ has label $\ob$.
Since $(1,5),(2,5),(3,5)\in D$ and $(1,7),(2,7),(3,7)\in D$, the box originating in $(4,5)$ must leave row $i$ before the box originating in $(6,7)$ enters row $i$ for $2\leq i \leq 4$. In particular, these boxes originating in $(4,5)$ and $(6,7)$ will only occupy the same row in row $1$.

By this fact, two of the boxes ${\sf b}_1,{\sf b}_2$ originating in$\{(4,4),(5,4),(6,4)\}$ must always lie weakly south of boxes originating in $(4,5)$ and $(6,7)$. Further, without loss of generality, ${\sf b}_1$ must always lie weakly south of the southmost $\gb$ in column $5$, and ${\sf b}_2$ must always lie weakly south of the southmost $\gb$ in column $7$.

Equation~\eqref{eq:fillBox} requires $\#(D\cap \{(4,5),(4,7)\})\geq 1.$ 
Thus, ${\sf b}_1$ must always lie weakly south of row $3$, and ${\sf b}_2$ must always lie weakly south of row $4$, or vice versa.
Since $(2,4)$ does not have label $\gb$, $\{(1,4),(2,4)\}\cap D\leq 1$, contradicting $\wt(D)=\gamma$. 
Then 
Equation~\ref{eq:failConj} follows. 
\end{proof}

By an exhaustive computer check, we have confirmed $w=12365847$ is a minimal counterexample to Conjecture~\ref{conj:main} with respect to $\ell(w)$ in $S_8$.
Our computations within $S_8$ suggest $\#\KKoh{w,\gamma}\leq g_{w,\gamma}$ in general.

Grothendieck polynomials for vexillary, \emph{i.e.}, $2143$-avoiding, permutations often have tamer combinatorial descriptions than the general case, see \cite{KMY,BPD}. In particular, as noted by T.~Matsumura and S.~Sugimoto \cite{MS20}, vexillary Grothendieck polynomials are Lascoux polynomials.
Thus in light of \cite{LascKKoh},
one might hope Conjecture~\ref{conj:main} still holds in the vexillary case, 
but this is false. For example, Conjecture~\ref{conj:main} fails for $w=12375846$. In Section~\ref{sec:321pf}, we prove Conjecture~\ref{conj:main} for $321$-avoiding permutations.

\section{An updated conjecture}\label{sec:newRule}
Computational evidence suggests the Ross--Yong rule weakly undercounts $g_{w,\gamma}$. Thus we seek a suitable relaxation. 
Let $D\subseteq[n]\times [n]$ be a diagram. Boxes in $D$ are labeled with either $\ob$ or $\gb$. 
Suppose $(i,j)\in D$ is rightmost in row $i\in[n]$ and has label $\gb$.  Let \[i'=\max\{r\in[i] \, | \, (r,j)\not\in D\}.\] Then if each $(r,j)\in D$ for $i'+1\leq r \leq i$ has label $\gb$, define the following:
\begin{itemize}
    \item[--] The \mydef{ghost Kohnert move} on $D$ at $(i,j)$ outputs $D'=D-\{(i,j)\}\cup \{(i',j)\}$. The new box $(i',j)$ has label $\gb$.
    \item[--] The \mydef{ghost K-Kohnert move} on $D$ at $(i,j)$ outputs $D'=D\cup \{(i',j)\}$. The new box $(i',j)$ has label $\gb$.
\end{itemize}

\begin{example}
Let $D$ be the diagram to the left. The middle diagram is the output of applying a ghost Kohnert move to $(5,1)$. The rightmost diagram is the output of applying a ghost K-Kohnert move to $(5,1)$.
\[\begin{picture}(250,60)
\put(0,40){\ytableausetup
{boxsize=1em}
{\begin{ytableau}
    \ob &  \ob  &  \ \\
    \ob &  \ & \ob \\
     \ob  &  \  &  \ \\
     \ &  \  &  \ \\
      \gb &  \  & \  \\
\end{ytableau}}}
\put(100,40){\ytableausetup
{boxsize=1em}
{\begin{ytableau}
 \ob &  \ob  &  \ \\
    \ob &  \ & \ob \\
     \ob  &  \  &  \ \\
     \gb&  \  &  \ \\
      \ &  \  & \  \\
\end{ytableau}}}
\put(200,40){\ytableausetup
{boxsize=1em}
{\begin{ytableau}
  \ob &  \ob  &  \ \\
    \ob &  \ & \ob \\
     \ob  &  \  &  \ \\
     \gb &  \  &  \ \\
      \gb &  \  & \  \\
\end{ytableau}}}
\end{picture}
\]
\end{example}

Let $\tKKoh{D}$ denote the set of diagrams obtainable through successive Kohnert, K-Kohnert, ghost Kohnert, and ghost K-Kohnert moves on $D$, where the initial diagram $D$ has all boxes labelled $\ob$. Take \[\tKKoh{w,\gamma}:=\{D\in \tKKoh{D(w)} \, | \, \wt(D)=\gamma\}.\]

\begin{example}\label{ex:newDiag}
Below is $\tKKoh{w,\gamma}$ for $w=12365847$ and $\gamma=(3,3,3,2)$, the counterexample given in Claim~\ref{claim:main}.
    \[
\begin{picture}(350,80)
\put(0,65){\ytableausetup
{boxsize=1em}
{\begin{ytableau}
 \ &  \ &  \ &  \ob &  \ob & \ &  \ob &  \ \\
   \ &  \ &  \ &  \ob &  \gb & \ &  \gb &  \ \\
    \ &  \ &  \ &  \gb &  \gb & \ &  \gb &  \ \\
     \ &  \ &  \ &  \ob &  \gb & \ &  \ &  \ \\
   \ &  \ &  \ &  \ &  \ & \ &  \ &  \ \\
    \ &  \ &  \ &  \ &  \ & \ &  \ &  \ \\
\end{ytableau}}}
\put(125,65){\ytableausetup
{boxsize=1em}
{\begin{ytableau}
 \ &  \ &  \ &  \ob &  \ob & \ &  \ob &  \ \\
   \ &  \ &  \ &  \gb &  \gb & \ &  \gb &  \ \\
    \ &  \ &  \ &  \ob &  \gb & \ &  \gb &  \ \\
     \ &  \ &  \ &  \ob &  \gb & \ &  \ &  \ \\
   \ &  \ &  \ &  \ &  \ & \ &  \ &  \ \\
    \ &  \ &  \ &  \ &  \ & \ &  \ &  \ \\
\end{ytableau}}}
\put(250,65){\ytableausetup
{boxsize=1em}
{\begin{ytableau}
 \ &  \ &  \ &  \ob &  \ob & \ &  \ob &  \ \\
   \ &  \ &  \ &  \gb &  \gb & \ &  \gb &  \ \\
    \ &  \ &  \ &  \ob &  \gb & \ &  \gb &  \ \\
     \ &  \ &  \ &  \ob &  \ & \ &  \gb &  \ \\
   \ &  \ &  \ &  \ &  \ & \ &  \ &  \ \\
    \ &  \ &  \ &  \ &  \ & \ &  \ &  \ \\
\end{ytableau}}}
\end{picture}
\]
The leftmost diagram is the element in $\tKKoh{w,\gamma}-\KKoh{w,\gamma}$. 
Below we provide the details with how to obtain this leftmost diagram.
For succinctness, we omit empty columns. Additionally, we omit some intermediate diagrams which use only ordinary Kohnert moves. 
%
\begin{equation*}
\begin{picture}(450,180)
\put(0,160){\ytableausetup
{boxsize=1em}
{\begin{ytableau}
   \ &  \ &   \  \\
     \ &  \ &   \  \\
      \ &  \ &   \ \\
     \ob &  \ob  &  \ \\
    \ob &  \ &   \ \\
    \ob &  \ &   \ob  \\
\end{ytableau}}}
 \put(58,132){$\rightarrow$}
\put(80,160){\ytableausetup
{boxsize=1em}
{\begin{ytableau}
 \ob &  \ &   \  \\
     \ &  \ &   \  \\
      \ &  \ &   \ \\
     \ob &  \ob  &  \ \\
    \ &  \ &   \ \\
    \ob &  \ &   \ob  \\
\end{ytableau}}}
 \put(138,132){$\rightarrow$}
\put(160,160){\ytableausetup
{boxsize=1em}
{\begin{ytableau}
 \ob &  \ &   \  \\
     \ &  \ &   \  \\
      \ &  \ob &   \ \\
     \ob &  \gb  &  \ \\
    \ &  \ &   \ \\
    \ob &  \ &   \ob  \\
\end{ytableau}}}
 \put(218,132){$\rightarrow$}
\put(240,160){\ytableausetup
{boxsize=1em}
{\begin{ytableau}
 \ob &  &   \  \\
     \ & \ob &   \  \\
      \ &  \ &   \ \\
     \ob &  \gb  &  \ \\
    \ &  \ &   \ \\
    \ob &  \ &   \ob  \\
\end{ytableau}}}
 \put(298,132){$\rightarrow$}
\put(320,160){\ytableausetup
{boxsize=1em}
{\begin{ytableau}
 \ob & \ob &   \  \\
     \ & \gb &   \  \\
      \ &  \ &   \ \\
     \ob &  \gb  &  \ \\
    \ &  \ &   \ \\
    \ob &  \ &   \ob  \\
\end{ytableau}}}
 \put(378,132){$\rightarrow$}
\put(400,160){\ytableausetup
{boxsize=1em}
{\begin{ytableau}
 \ob & \ob &   \  \\
     \ & \gb &   \  \\
      \ &  \ &   \ \\
     \ob &  \gb  &   \ob \\
    \ &  \ &   \ \\
    \ob &  \ &  \  \\
\end{ytableau}}}
\put(420,80){$\downarrow$}
\put(400,50){\ytableausetup
{boxsize=1em}
{\begin{ytableau}
\ob & \ob &   \  \\
     \ & \gb &   \  \\
      \ob &  \ &   \ \\
     \ob &  \gb  &   \ob \\
    \ &  \ &   \ \\
    \ &  \ &  \  \\
\end{ytableau}}}
 \put(378,22){$\leftarrow$}
\put(320,50){\ytableausetup
{boxsize=1em}
{\begin{ytableau}
\ob & \ob &   \  \\
     \ob & \gb &   \  \\
      \gb &  \ &   \ \\
     \ob &  \gb  &   \ob \\
    \ &  \ &   \ \\
    \ &  \ &  \  \\
\end{ytableau}}}
 \put(298,22){$\leftarrow$}
\put(240,50){\ytableausetup
{boxsize=1em}
{\begin{ytableau}
\ob & \ob &   \  \\
     \ob & \gb &   \  \\
      \gb &  \ &   \ob \\
     \ob &  \gb  &   \ \\
    \ &  \ &   \ \\
    \ &  \ &  \  \\
\end{ytableau}}}
 \put(218,22){$\leftarrow$}
\put(160,50){\ytableausetup
{boxsize=1em}
{\begin{ytableau}
\ob & \ob &   \  \\
     \ob & \gb &   \  \\
      \gb &  \gb &   \ob \\
     \ob &  \gb  &   \ \\
    \ &  \ &   \ \\
    \ &  \ &  \  \\
\end{ytableau}}}
 \put(138,22){$\leftarrow$}
\put(80,50){\ytableausetup
{boxsize=1em}
{\begin{ytableau}
\ob & \ob &   \  \\
     \ob & \gb &   \ob  \\
      \gb &  \gb &   \gb \\
     \ob &  \gb  &   \ \\
    \ &  \ &   \ \\
    \ &  \ &  \  \\
\end{ytableau}}}
 \put(58,22){$\leftarrow$}
\put(00,50){\ytableausetup
{boxsize=1em}
{\begin{ytableau}
\ob & \ob &   \ob  \\
     \ob & \gb &   \gb  \\
      \gb &  \gb &   \gb \\
     \ob &  \gb  &   \ \\
    \ &  \ &   \ \\
    \ &  \ &  \  \\
\end{ytableau}}}
\end{picture}
\end{equation*}
\end{example}

\begin{conjecture}\label{conj:newKKohGr}
   For $w\in S_n$ and $\gamma\in\mathbb{Z}_{\geq 0}^n$,  $g_{w,\gamma}=\#\tKKoh{w,\gamma}$. Thus,
   \[{\mathfrak{G}_w}= \sum_{D\in \tKKoh{D(w)}}(-1)^{\#D-\ell(w)}x^{\wt(D)}.\]
\end{conjecture}

 Conjecture~\ref{conj:newKKohGr} has been checked exhaustively for $S_n$ with $n\leq 8$.  In the next section we prove Conjecture~\ref{conj:newKKohGr} for $321$-avoiding permutations. Intuitively, Conjecture~\ref{conj:newKKohGr} seems promising, as the additional ghost moves allow this rule to more closely mirror the established pipe dream formula for Grothendieck polynomials \cite{FominKrillov}.

J.~Pan and T.~Yu \cite{LascKKoh} prove Conjecture~\ref{thm:Lasc}, \emph{i.e.}, that $\KKoh{D(\alpha)}$ generate ${\mathfrak{L}_\alpha}$. 
After private communication, T.~Yu proves that the same holds for $\tKKoh{D(\alpha)}$, using a straightforward analysis of \cite[Theorem~44]{LascKKoh}:
\begin{proposition}\label{prop:newKKohLasc}\cite{YuLasc}
  For $\alpha\in\mathbb{Z}_{\geq 0}^n$, $\KKoh{D(\alpha)}=\tKKoh{D(\alpha)}$.
Thus $\tKKoh{D(\alpha)}$ generates Lascoux polynomial ${\mathfrak{L}_\alpha}$. 
\end{proposition}

In fact, suppose we further relax Conjecture~\ref{conj:newKKohGr} so that one may apply ghost Kohnert and ghost K-Kohnert moves to 
\emph{any} $(i,j)\in D$ with label $\gb$. Take $\bKKoh{D}$  to be the corresponding set of diagrams generated by these relaxed moves, along with the original Kohnert and K-Kohnert moves. T.~Yu proves $\bKKoh{D(\alpha)}$ generates  ${\mathfrak{L}_\alpha}$.
Computational evidence suggests $\bKKoh{D(w)}$ weakly overcounts $g_{w,\gamma}$. 
For example,  $w=12385746$ and $\gamma=(1,4,1,2)$, there are $8$ diagrams in $\bKKoh{D(w)}$ with weight $\gamma$, but  $g_{w,\gamma}=7$. However for $w\in S_7$, $\bKKoh{D(w)}$ correctly computes $g_{w,\gamma}$.

\section{Proofs for $321$-avoiding permutations}\label{sec:321pf}
In this section we prove these conjectures in the $321$-avoiding case. We utilize a correspondence between diagrams generating Grothendieck polynomials and Lascoux polynomials.
\begin{theorem}\label{thm:old321} Conjecture~\ref{conj:main} and Conjecture~\ref{conj:newKKohGr} hold for $w\in S_n$ $321$-avoiding.
\end{theorem}
We first provide additional combinatorial background. Then we end with the proof of Theorem~\ref{thm:old321}.

\subsection{Flagged set-valued tableaux}
A \mydef{flagged set-valued tableau} for $D\subseteq[n]\times[n]$ is a filling $f:D\rightarrow 2^{[n]}-\emptyset$ such that 
    \begin{itemize}
        \item $\min f(r,c)\geq \max f(r,c+k)$ for $(r,c),(r,c+k)\in D$ where $k>0$.
        \item $\max f(r,c)< \min f(r+k,c)$ for $(r,c),(r+k,c)\in D$ where $k>0$.
        \item $\max f(r,c)\leq r$ for $(r,c)\in D$.
    \end{itemize}
Here $2^{[n]}$ indicates the power set of $[n]$.    
Let ${\sf FSVT}(D)$ denote the set of flagged set-valued tableau for $D$. 
The \mydef{weight} of $T$ in ${\sf FSVT}(D)$, denoted $\wt(T)\in \mathbb{Z}_{\geq 0}^n$, is defined by
$\wt(T)_i:=\# i\text{'s} \text{ in } T$.

    T.~Matsumura \cite{Mats20} gives a formula for $321$-avoiding Grothendieck polynomials in terms of flagged set-valued tableaux:
   \begin{theorem}\label{thm:mats321}\cite[Theorem 3.1]{Mats20} For $w\in S_n$ $321$-avoiding,
     \[{\mathfrak{G}_w}= \sum_{T\in {\sf{FSVT}}(D(w))}(-1)^{\#T-\ell(w)}x^{\wt(T)},\]
     where $\#T=\sum_{i\in[n]} \wt(T)_i$.
\end{theorem}
For the remainder of this section, assume $w\in S_n$ is $321$-avoiding. 
Let 
\[\overline{D(w)}:=\{(i,j) \ | \ (i,j+k),(i-k',j)\in D(w) \text{ for some } k,k'\geq 0\},\]
\emph{i.e.}, the southwest hook closure of $D(w)$.
Define the map $\phi:{\sf{FSVT}}(D(w))\rightarrow {\sf{FSVT}}(\overline{D(w)})$, where for $T\in {\sf{FSVT}}(D(w))$ and $(r,c)\in \overline{D(w)}$,
\begin{align*}
    (\phi(T))(r,c):=
    \begin{cases}
        T(r,c) & \text{ if } (r,c)\in D(w)\\
        \{r\} & \text{ else.} 
    \end{cases}
\end{align*}

\begin{claim}
  The map $\phi$ is a well-defined injection.   
\end{claim}
\begin{proof}
    Consider $(i,j)\in \overline{D(w)}-D(w)$. Since $(i,j)\in \overline{D(w)}$, we know $w(i)>j$ and there exists some $i'<i$ such that $w(i')>j$.
Further, since $(i,j)\not\in D(w)$, there exists $i'<i''<i$ such that $w(i'')=j$. 

Suppose there exists $(i,j')\in D(w)$ for some $j'<j$. By the definition of $D(w)$, $w(h)=j'$ for some $h>i$. Then $i'<i''<h$ and $w(i')>w(i'')>w(h)$, giving a $321$-pattern in $w$. 

Thus for $w$ is $321$-avoiding,  $(i,j)\in \overline{D(w)}-D(w)$ implies that if $j'<j$, $(i,j')\not\in D(w)$.  Then by construction, $\phi(T)\in {\sf{FSVT}}(\overline{D(w)})$. Injectivity follows by the construction of $\phi$.
\end{proof}

Define $\al{w},\be{w}\in\mathbb{Z}_{\geq 0}^n$ such that
\begin{align*}
    (\al{w})_i:=&\#\{(i,j)\in \overline{D(w)}\}\\
    (\be{w})_i:=&\#\{(i,j)\in {D(w)}\}.
\end{align*}


By construction,   
\[(\al{w})_i=\#\{(i,j) \, | \, (i',j)\in D(w) \text{ for some } i'\leq i\}.\] Thus the nonzero parts of $\al{w}$ are weakly increasing. 
Similarly, \[(\al{w})_i-(\be{w})_i=\#\{(i,j) \, | \, (i',j)\in D(w) \text{ for some } i'< i, \text{ where } (i,j)\not\in D(w)\},\] so
the nonzero parts $\al{w}-\be{w}$ are also weakly increasing. Note for $T\in {\sf{FSVT}}(D(w))$, $\wt(T)=\wt(\phi(T))-(\al{w}-\be{w})$.

Let $\psi_w:\KKoh{\overline{D(w)}}\rightarrow \KKoh{D(\al{w})}$ be the map that deletes empty columns in the grid. 
Restricting $\psi_w$ to $\overline{D(w)}$ induces a bijection $\rho_w:{\sf{FSVT}}(\overline{D(w)})\rightarrow{\sf{FSVT}}(D(\al{w}))$.
\begin{example}\label{ex:FSVText}
In each of the tableaux below, we have shaded the underlying diagram.
To the left is some $T\in {\sf{FSVT}}(D(w))$ for $w=451829367$. The middle tableau is $\phi(T)$. To the right is $\rho_w(\phi(T))$:
\[\begin{picture}(380,75)
\put(0,60){\ytableausetup
{boxsize=1.2em}
{\begin{ytableau}
 *(gray!50) 1 &  *(gray!50)1 &  *(gray!50) 1&  \ &  \ & \ &  \ \\
 *(gray!50) 2&  *(gray!50)2 &  *(gray!50) 2&  \ &  \ & \ &  \ \\
 \ &  \ &  \ &  \ &  \ & \  &  \  \\
 \ &  *(gray!50) 4  &  *(gray!50) 4 &  \ &  \ &  *(gray!50) \scriptstyle{4 3 2}&  *(gray!50) 21\\
\ &  \ &  \ &  \ &  \ & \  &  \  \\
\ &  \ &  *(gray!50)65 &  \ &  \ &  *(gray!50)5 &  *(gray!50)\scriptstyle{543}  \\
\end{ytableau}}}
\put(150,60){\ytableausetup
{boxsize=1.2em}
{\begin{ytableau}
 *(gray!50) 1 &  *(gray!50)1 &  *(gray!50) 1&  \ &  \ & \ &  \ \\
 *(gray!50) 2&  *(gray!50)2 &  *(gray!50) 2&  \ &  \ & \ &  \ \\
 \ &  \ &  \ &  \ &  \ & \  &  \  \\
 *(gray!50)4 &  *(gray!50) 4  &  *(gray!50) 4 &  \ &  \ &  *(gray!50) \scriptstyle{4 3 2}&  *(gray!50) 21\\
\ &  \ &  \ &  \ &  \ & \  &  \  \\
*(gray!50)6 &  *(gray!50)6 &  *(gray!50)65 &  \ &  \ &  *(gray!50)5 &  *(gray!50)\scriptstyle{543}  \\
\end{ytableau}}}
\put(300,60){\ytableausetup
{boxsize=1.2em}
{\begin{ytableau}
 *(gray!50) 1 &  *(gray!50)1 &  *(gray!50) 1&  \ &  \  \\
 *(gray!50) 2&  *(gray!50)2 &  *(gray!50) 2&  \ &  \  \\
 \ &  \ &  \ &  \ &  \ \\
 *(gray!50)4 &  *(gray!50) 4  &  *(gray!50) 4 &  *(gray!50) \scriptstyle{4 3 2}&  *(gray!50) 21  \\
\ &  \ &  \ &  \ &  \  \\
*(gray!50)6 &  *(gray!50)6 &  *(gray!50)65  &  *(gray!50)5 &  *(gray!50)\scriptstyle{543}  \\
\end{ytableau}} \ \raisebox{-2.5cm}{.}}
\end{picture}
\]
\end{example}

\subsection{Set-valued key tableaux}
As proven in \cite{LascKKoh}, $\KKoh{D(\alpha)}$ bijects with particular set-valued tableaux. We describe this correspondence, specialized to the case in which the nonzero parts in $\alpha$ are weakly increasing. Assume $D(\alpha)\subseteq [n]\times[n]$.

A \mydef{tableau} for $D(\alpha)$ is a filling $f:D(\alpha)\rightarrow [n]$ such that 
    \begin{itemize}
        \item $f(r,c)\geq f(r,c+k)$ for $(r,c),(r,c+k)\in D(\alpha)$ where $k>0$.
        \item $f(r,c)< f(r+k,c)$ for $(r,c),(r+k,c)\in D(\alpha)$ where $k>0$.
    \end{itemize}

Let ${\sf Tab}(D(\alpha))$ denote the set of tableaux for $D(\alpha)$.
For $T\in {\sf FSVT}(D(\alpha))$, define the tableau $M(T)\in {\sf Tab}(D(\alpha))$ such that \[(M(T))(r,c)=\max T(r,c) \text{ for } (r,c)\in D(\alpha).\]
The \mydef{left key} of $T$, denoted $K_{-}(T)\in {\sf Tab}(D(\alpha))$, is defined recursively. 
We construct the columns $\{C_k\}_{k\in[n]}$ of $K_{-}(T)$, left to right, using the algorithm of \cite{Willis}.

Take $T_1(k)$ be $M(T)$ restricted to its $k$ leftmost columns. 
Let 
\[m_k=\# \text{ boxes in column } k \text{ of }M(T).\]
For $j\in[m_k]$ in increasing order, we compute the sequence $c^{(j)}=(c_i^{(j)})_{i\in[k]}$ from $T_j(k)$.
Define $c_1^{(j)}$ as the southmost entry in column $k$ of $T_j(k)$. Take $c_{i+1}^{(j)}$ to be the smallest entry in column $k+1-i$ of $T_j(k)$ that is weakly greater than $c_i^{(j)}$, where $i\in[k-1]$. Append the entry $c_k^{(j)}$ to the set of column entries $C_k$. Then:
\begin{itemize}
    \item If $j<m_k$, let $T_{j+1}(k)$ be the tableau formed from $T_j(k)$ by deleting any entries in column $i$ weakly south of $c_i^{(j)}$ for each $i\in[k]$. 
    \item If $j=m_k$, place entries of $C_k$ into the $k$th column of $K_{-}(T)$ such that entries increase down columns.
\end{itemize}

\begin{example}
To the left is some $T\in {\sf FSVT}(D(\alpha))$ for $\alpha=(2,3,0,4,0,4)$ and to the right is the corresponding $K_{-}(T)$. The empty rows in $D(\alpha)$ have been deleted throughout.
 \[
\begin{picture}(225,80)
\put(0,55){\ytableausetup{boxsize=0.52cm}
\begin{ytableau}
*(gray!50) 1 &  *(gray!50)1 &  \ &  \ \\
 *(gray!50) 2&  *(gray!50)2 &  *(gray!50) 1&  \ \\
 *(gray!50) 43 &  *(gray!50) 3  &  *(gray!50) 3&    *(gray!50) 32\\
*(gray!50)6 &  *(gray!50)6 &  *(gray!50)65 &   *(gray!50)54   \\
\end{ytableau}}
\put(150,55){\ytableausetup{boxsize=0.52cm}
\begin{ytableau}
*(gray!50) 1 &  *(gray!50)1 &  \ &  \ \\
 *(gray!50) 2&  *(gray!50)2 &  *(gray!50) 1&  \ \\
 *(gray!50) 4 &  *(gray!50) 4 &  *(gray!50) 4 &    *(gray!50) 4 \\
*(gray!50)6 &  *(gray!50)6 &  *(gray!50)6 &   *(gray!50)6   \\
\end{ytableau}}
\end{picture}\]
Below are the tableaux $T_j(k)$ with the corresponding $k$ increasing left to right.  In each, the sequences ${c^{(j)}}$ are noted in blue with the particular entries $c_k^{(j)}$ in bold.
\[
\begin{picture}(300,200)
\put(0,160){$T_1(k)$:}
\put(50,180){
\ytableausetup{boxsize=0.45cm}
\begin{ytableau}
 *(gray!50) 1 \\
 *(gray!50) 2 \\
 *(gray!50) 4 \\
*(gray!50){\textcolor{blue}{\mathbf{6}}}  \\
\end{ytableau}}
\put(100,180){\ytableausetup{boxsize=0.45cm}
\begin{ytableau}
*(gray!50) 1 &  *(gray!50)1\\
 *(gray!50) 2&  *(gray!50)2 \\
 *(gray!50) 4 &  *(gray!50) 3   \\
*(gray!50){\textcolor{blue}{\mathbf{6}}} &  *(gray!50){\textcolor{blue}{6}}  \\
\end{ytableau}\
}
\put(160,180){\ytableausetup{boxsize=0.45cm}
\begin{ytableau}
*(gray!50) 1 &  *(gray!50)1  &  \  \\
 *(gray!50) 2&  *(gray!50)2 &  *(gray!50) 1  \\
 *(gray!50) 4  &  *(gray!50) 3   &  *(gray!50) 3  \\
*(gray!50){\textcolor{blue}{\mathbf{6}}} &  *(gray!50){\textcolor{blue}{6}} &  *(gray!50){\textcolor{blue}{6}}    \\
\end{ytableau}}
\put(240,180){\ytableausetup{boxsize=0.45cm}
\begin{ytableau}
*(gray!50) 1 &  *(gray!50)1 &  \ &  \ \\
 *(gray!50) 2&  *(gray!50)2 &  *(gray!50) 1&  \ \\
 *(gray!50) 4  &  *(gray!50) 3  &  *(gray!50) 3  &    *(gray!50) 3  \\
*(gray!50){\textcolor{blue}{\mathbf{6}}} &  *(gray!50){\textcolor{blue}{6}} &  *(gray!50){\textcolor{blue}{6}} &   *(gray!50){\textcolor{blue}{5}}   \\
\end{ytableau}}
\put(0,95){$T_2(k)$:}
\put(50,105){
\ytableausetup{boxsize=0.45cm}
\begin{ytableau}
 *(gray!50) 1  \\
 *(gray!50) 2 \\
 *(gray!50) {\textcolor{blue}{\mathbf{4}}} \\
\end{ytableau}}
\put(100,105){\ytableausetup{boxsize=0.45cm}
\begin{ytableau}
*(gray!50) 1  &  *(gray!50)1 \\
 *(gray!50) 2 &  *(gray!50)2  \\
 *(gray!50) {\textcolor{blue}{\mathbf{4}}} &  *(gray!50) {\textcolor{blue}{3}}  \\
\end{ytableau}\
}
\put(160,105){\ytableausetup{boxsize=0.45cm}
\begin{ytableau}
*(gray!50) 1  &  *(gray!50)1  &  \  \\
 *(gray!50) 2&  *(gray!50)2 &  *(gray!50) 1 \\
 *(gray!50) {\textcolor{blue}{\mathbf{4}}} &  *(gray!50) {\textcolor{blue}{3}}  &  *(gray!50) {\textcolor{blue}{3}}\\
\end{ytableau}}
\put(240,105){\ytableausetup{boxsize=0.45cm}
\begin{ytableau}
*(gray!50) 1 &  *(gray!50)1 &  \ &  \ \\
 *(gray!50) 2&  *(gray!50)2 &  *(gray!50) 1&  \ \\
 *(gray!50) {\textcolor{blue}{\mathbf{4}}} &  *(gray!50) {\textcolor{blue}{3}}  &  *(gray!50) {\textcolor{blue}{3}} &    *(gray!50) {\textcolor{blue}{3}} \\
\end{ytableau}}
\put(0,40){$T_3(k)$:}
\put(50,45){
\ytableausetup{boxsize=0.45cm}
\begin{ytableau}
 *(gray!50) 1 \\
 *(gray!50) {\textcolor{blue}{\mathbf{2}}} \\
\end{ytableau}}
\put(100,45){\ytableausetup{boxsize=0.45cm}
\begin{ytableau}
*(gray!50) 1 &  *(gray!50)1\\
 *(gray!50) {\textcolor{blue}{\mathbf{2}}}&  *(gray!50){\textcolor{blue}{2}} \\
\end{ytableau}\
}
\put(160,45){\ytableausetup{boxsize=0.45cm}
\begin{ytableau}
*(gray!50) {\textcolor{blue}{\mathbf{1}}} &  *(gray!50){\textcolor{blue}{1}} &  \  \\
 *(gray!50) 2&  *(gray!50)2 &  *(gray!50) {\textcolor{blue}{1}} \\
\end{ytableau}}
\put(0,5){$T_4(k)$:}
\put(50,0){
\ytableausetup{boxsize=0.45cm}
\begin{ytableau}
 *(gray!50) {\textcolor{blue}{\mathbf{1}}} \\
\end{ytableau}}
\put(100,0){\ytableausetup{boxsize=0.45cm}
\begin{ytableau}
*(gray!50) {\textcolor{blue}{\mathbf{1}}} &  *(gray!50){\textcolor{blue}{1}}\\
\end{ytableau}\
}
\end{picture}
\]
\end{example}

Take $\alpha\in\mathbb{Z}_{\geq 0}^n$ whose nonzero parts are weakly increasing.
 A \mydef{set-valued key tableau} for $\alpha$ 
 is a filling $f:D(\alpha)\rightarrow 2^{[n]}$ such that 
    \begin{enumerate}
        \item $\min f(r,c)\geq \max f(r,c+k)$ for $(r,c),(r,c+k)\in D(\alpha)$ where $k>0$.
        \item $\max f(r,c)< \min f(r+k,c)$ for $(r,c),(r+k,c)\in D(\alpha)$ where $k>0$.
        \item $K_{-}(T)(r,c)\leq r$ for $(r,c)\in D(\alpha)$.
    \end{enumerate}
 
  Let ${\sf SVKT}(D(\alpha))$ denote the set of set-valued key tableau for $\alpha$.
     By the definition of the left key, it is immediate that ${\sf{FSVT}}(D(\alpha))={\sf{SVKT}}(D(\alpha))$. 
This produces the following specialization:
\begin{theorem}\cite{BSW,LascKKoh,SU}\label{thm:lascExp}
     For a weak composition $\alpha\in\mathbb{Z}_{\geq 0}^n$ such that the nonzero parts of $\alpha$ are weakly increasing,
     \[{\mathfrak{L}_\alpha}=\hspace{-.15cm} \sum_{T\in {\sf{SVKT}}(D(\alpha))}\hspace{-.1cm} (-1)^{\#T-|\alpha|}x^{\wt(T)}=\hspace{-.15cm} \sum_{T\in {\sf{FSVT}}(D(\alpha))}\hspace{-.1cm}(-1)^{\#T-|\alpha|}x^{\wt(T)}=\hspace{-.15cm}\sum_{T\in \KKoh{D(\alpha})}\hspace{-.1cm}(-1)^{\#D-|\alpha|}x^{\wt(D)},\]
     where $|\alpha|=\sum_{i\in[n]}\alpha_i$.
\end{theorem}

Take $\alpha\in\mathbb{Z}_{\geq 0}^n$ such that the nonzero parts of $\alpha$ are weakly increasing.
We describe the weight-preserving bijection $\Phi_{\alpha}:{\sf{FSVT}}(D(\alpha))\rightarrow \KKoh{D(\alpha)}$ in  \cite{LascKKoh}
that gives the last equality in Theorem~\ref{thm:lascExp}.

Encode $T\in {\sf{FSVT}}(D(\alpha))$ as a pair of diagrams $(O(T),G(T))$ in $[n]\times [n]$, where 
\begin{align}\label{eq:tabPair}
\begin{split}
    O(T)&=\{(i,j) \ | \ i=\max T(r,j) \text{ for some } r\in[n]\}, \text{ and}\\
    G(T)&=\{(i,j) \ | \ i\in T(r,j)-\{\max T(r,j)\} \text{ for some } r\in[n]\}.
\end{split}
\end{align}

Define the map $\Phi_{\alpha}$ by the following algorithm:
\begin{itemize}
    \item Suppose $T\in {\sf{FSVT}}(D(\alpha))$. Initialize $S:=O(T)$.
    Iterate over $G(T)$ upwards in columns, working from left to right. 
    \item For each $(i,j)\in G(T)$, pick the minimal $i'\geq i$ such that $(i',j)\in S$. Update $S$ to be $S-\{(i',j)\}\cup\{(i,j)\}$.
    \item After iterating through $G(T)$, output the labelled diagram $\Phi_{\alpha}(T)=O'\cup G'$, where
    $O':=S$ and $G':=O(T)\cup G(T)-S$. Here, boxes in $O'$ are assigned label $\ob$, and boxes in $G'$ are assigned label $\gb$. 
\end{itemize}

We describe the inverse map $\Phi_{\alpha}^{-1}$ for completeness:
\begin{itemize}
    \item Suppose $D\in \KKoh{D(\alpha)}$. Let $O(D)$ be the boxes in $D$ with label $\ob$ and $G(D)$ be the boxes in $D$ with label $\gb$.
    \item Initialize $S:=O(D)$.
    Iterate through boxes in $G(D)$ down columns, working from right to left. 
    \item For each $(i,j)\in G(D)$, pick the maximal $i'\leq i$ such that $(i',j)\in S$. Update $S$ to be $S-\{(i',j)\}\cup\{(i,j)\}$.
    \item After iterating through $G(D)$, output the (encoded) tableau $\Phi_{\alpha}(D)=(O',G')$, where
    $O':=S$ and $G':=O(D)\cup G(D)-S$.
\end{itemize}

\begin{example}\label{ex:bij}
The leftmost diagram is the pair $(O(U),G(U))$, where $U=\rho_w(\phi(T))$ as in Example~\ref{ex:FSVText}. To the right is $\Phi_{\alpha}(U)$.
    \[\begin{picture}(250,75)
\put(0,60){\ytableausetup
{boxsize=1.2em}
{\begin{ytableau}
 \ob &  \ob &  \ob      &  \ &  \gb  \\
 \ob &  \ob &  \ob      &  \gb &  \ob  \\
 \ & \ &        \  &      \gb    &   \gb \\
   \ob &  \ob &  \ob    &  \ob &  \gb    \\ 
   \ & \  &  \gb &  \ob   &  \ob         \\ 
   \ob &  \ob &  \ob    &  \ &  \  \\
\end{ytableau}}}
\put(150,60){\ytableausetup
{boxsize=1.2em}
{\begin{ytableau}
 \ob    &  \ob  &  \ob  &  \ & \ob   \\
 \ob    &  \ob  &  \ob  &\ob &  \gb \\
 \      & \     & \     & \gb&   \ob \\
   \ob  &  \ob  &  \ob  &\gb &  \gb  \\ 
   \    & \     &  \ob  &\ob &  \gb \\ 
   \ob  &  \ob  &  \gb  &  \ &  \    \\
\end{ytableau}}}
\end{picture}.
\]
\end{example}

\subsection{Proof of Theorem~\ref{thm:old321}}\label{sec:mainPf}
Assume $w\in S_n$ is $321$-avoiding. 
We will construct a weight-preserving bijection $f:{\sf{FSVT}}(D(w))\rightarrow \KKoh{D(w)}$.

First consider the weight-preserving bijection $\rho_w:S_1\rightarrow S_2$, where
\begin{align*}
    S_1:=&\{\phi(T)  \ | \ T\in {\sf{FSVT}}(D(w))\}, \text{ and}\\
    S_2:=&\{T\in {\sf{FSVT}}(D(\al{w})) \ | \ T(i,j)=i \text{ if } j\leq (\al{w})_i-(\be{w})_i\}.
\end{align*}
Since the nonzero parts of $\al{w}$ and $\al{w}-\be{w}$ are weakly increasing, the boxes 
\[\{(i,j)\in D(\al{w}) \ | \ j\leq (\al{w})_i-(\be{w})_i  \}\]
are maximally southwest in $D(\al{w})$.

Then consider $U\in S_2$ encoded as $(O(U),G(U))$ using Equation~\eqref{eq:tabPair}.
Suppose $i\in[n]$ is such that if  $(i',j)\in D(\al{w})$ with $i'\geq i$, then $\#U(i',j)=1$. Thus $(U(i',j),j)\in O(U)$. Now suppose $(r,j)\in G(U)$ for some $r\in[n]$. Then $(r,j)\in U(k,j)$ for some $k\in[n]$, where $(\max U(k,j),j)\in O(U)$. By assumption, $r<\max U(k,j)<i$.
Thus, $(U(i,j),j)$ will never be removed from $S$ when computing $\Phi_{\al{w}}$, so $(U(i,j),j)$ has label $\ob$ in $\Phi_{\al{w}}(U)$. 
This implies $\Phi_{\al{w}}(S_2)\subseteq S_3$, where
\[S_3:=\{D\in \KKoh{D(\al{w})} \ | \ (i,j)\in D \text{ with label } \ob \text{ for } j\leq (\al{w})_i-(\be{w})_i \}.\]

Now take $D\in S_3$ where $O(D)$ are the boxes labelled $\ob$ and $G(D)$ are those labelled $\gb$. Then consider $\Phi_{\al{w}}^{-1}(D)$.
Suppose $i\in[n]$ is such that for $i'\geq i$, if $(i',j)\in D$, then $(i',j)\in O(D)$. By the definition of $\Phi_{\al{w}}^{-1}$, a box $(r,j)$ might be removed from $S$ only if there exists some $(r',j)\in G(D)$ such that $r'>r$. This ensures $(i',j)\in O'$ where $\Phi_{\al{w}}^{-1}(D)=(O',G')$. 
Thus $\Phi_{\al{w}}(S_2)=S_3$.

Recall $\psi_w:\KKoh{\overline{D(w)}}\rightarrow \KKoh{D(\al{w})}$ is the map that deletes empty columns in the grid.
Take the map $\Psi:S_3\rightarrow \KKoh{D(w)}$ where
\[\Psi(D)=\psi_w^{-1}(D)-\{(i,j)\in \psi_w^{-1}(D) \ | \ (i,j)\in \overline{D(w)}-D(w)\}.\]
Note $\Psi$ is a bijection such that $\wt(D)=\wt(\Psi(D))+(\al{w}-\be{w})$.

Therefore, $f:=\Psi\circ\Phi_{\al{w}}\circ\rho_w\circ\phi$ is as desired, so Conjecture~\ref{conj:main} follows by Theorem~\ref{thm:mats321}.
By Proposition~\ref{prop:newKKohLasc}, Conjecture~\ref{conj:newKKohGr} follows by the same argument.  \qed

\begin{example}\label{ex:backToKoh}
   Take $U'=\Phi_{\alpha}(U)$ as given in Example~\ref{ex:bij}. On the right is $D=\Psi(U')$, where $\Psi$ is as defined in Section~\ref{sec:mainPf}. 
    \[\begin{picture}(275,95)
\put(0,75){\ytableausetup
{boxsize=1.2em}
{\begin{ytableau}
 *(gray!50) 1 &  *(gray!50)1 &  *(gray!50) 1&  \ &  \ & \ &  \ \\
 *(gray!50) 2&  *(gray!50)2 &  *(gray!50) 2&  \ &  \ & \ &  \ \\
 \ &  \ &  \ &  \ &  \ & \  &  \  \\
 \ &  *(gray!50) 4  &  *(gray!50) 4 &  \ &  \ &  *(gray!50) \scriptstyle{4 3 2}&  *(gray!50) 21\\
\ &  \ &  \ &  \ &  \ & \  &  \  \\
\ &  \ &  *(gray!50)65 &  \ &  \ &  *(gray!50)5 &  *(gray!50)\scriptstyle{543}  \\
\end{ytableau}}}
\put(125,35){\Large$\xrightarrow{ f}$}
\put(150,75){\ytableausetup
{boxsize=1.2em}
{\begin{ytableau}
 \ob    &  \ob  &  \ob &  \ &  \  &  \ & \ob   \\
 \ob    &  \ob  &  \ob &  \ &  \  &\ob &  \gb  \\
 \      & \     & \    &  \ &  \  & \gb&   \ob \\
   \ &  \ob  &  \ob &  \ &  \  &\gb &  \gb   \\ 
   \    & \     &  \ob &  \ &  \  &\ob &  \gb  \\ 
     &    &  \gb &  \ &  \  &  \ &  \    \\
\end{ytableau}}}
\end{picture}
\]
For $T$ as in Example~\ref{ex:FSVText}, recalled above, we find $f(T)=(\Psi\circ\Phi_{\al{w}}\circ\rho_w\circ\phi)(T)=D$. 
\end{example}

\begin{remark}
Recall from Section~\ref{sec:newRule} that $\bKKoh{D}$ is the set of diagrams generated by the original Kohnert and K-Kohnert moves while also allowing ghost moves on arbitrary boxes with label $\gb$.  
    By the argument in Theorem~\ref{thm:old321} along with \cite{YuLasc}, $\bKKoh{D(w)}$ generates $\mathfrak{G}_w$ for $321$-avoiding $w\in S_n$.
\end{remark}

 We expect Conjecture~\ref{conj:newKKohGr} may have tableaux-based proofs in the $1432$-avoiding and vexillary cases.
 In particular, one might mimic the J.~Pan and T.~Yu \cite{LascKKoh} set-valued tableaux argument using the rules of A.~Knutson, E.~Miller, and A.~Yong \cite{KMY} as well as N.~J.~Y.~Fan and P.~L.~Guo \cite{Fan.Guo} to prove Conjecture~\ref{conj:newKKohGr} for vexillary and $1432$-avoiding permutations, respectively.

\section*{Acknowledgements}
The author would like to thank Igor Pak, Avery St.~Dizier, Anna Weigandt, Alexander Yong, and Tianyi Yu for helpful comments and conversations. We also wish to thank Jianping Pan and Tianyi Yu for sharing their code to generate Kohnert diagrams. We also thank the referees for their thoughtful feedback.

\bibliographystyle{plainurl}
\bibliography{KKohnert}

\end{document}